\documentclass[a4paper,12pt]{article}
\usepackage{amsmath}
\usepackage{amsfonts}
\usepackage{amssymb}
\usepackage{latexsym}
\usepackage{epsfig}
\usepackage{graphicx}
\usepackage{oldgerm}
\usepackage{theorem}

\newcommand{\mib}[1]{\mbox{\boldmath $#1$}}
\def\R{\mathbb{R}}
\def\W{\mathbb{W}}
\def\Z{\mathbb{Z}}
\def\C{\mathbb{C}}
\def\N{\mathbb{N}}

\def\P{{\bf P}}
\def\E{{\bf E}}

\def\x{\mib{x}}
\def\y{\mib{y}}
\def\a{{\bf a}}
\def\b{{\bf b}}
\def\0{{\bf 0}}
\def\n{{\bf n}}
\def\X{\mib{X}}
\def\Y{\mib{Y}}

\def\cN{{\cal N}}

\def\Pf{\mathop{\mathrm{Pf}}}
\def\={\stackrel{\rm d}{=}}
\theorembodyfont{\itshape}

\newtheorem{thm}{Theorem}[section]
\newtheorem{lem}[thm]{Lemma}

\newcommand{\SSC}[1]{\section{#1}\setcounter{equation}{0}}

\begin{document}

\title{\bf Extreme value distributions \\
of noncolliding diffusion processes}
\author{
Minami IZUMI
\footnote{
Department of Physics,
Faculty of Science and Engineering,
Chuo University, 
Kasuga, Bunkyo-ku, Tokyo 112-8551, Japan
}
and
Makoto KATORI
\footnote{
Department of Physics,
Faculty of Science and Engineering,
Chuo University, 
Kasuga, Bunkyo-ku, Tokyo 112-8551, Japan;
e-mail: katori@phys.chuo-u.ac.jp
}}
\date{7 July 2010}
\pagestyle{plain}
\maketitle
\begin{abstract}
Noncolliding diffusion processes reported in the 
present paper are $N$-particle systems of diffusion processes
in one-dimension, which are conditioned so that
all particles start from the origin and
never collide with each other in a finite time interval $(0, T)$,
$0 < T < \infty$.
We consider four temporally inhomogeneous
processes with duration $T$,
the noncolliding Brownian bridge,
the noncolliding Brownian motion,
the noncolliding three-dimensional Bessel bridge,
and the noncolliding Brownian meander.
Their particle distributions at each time $t \in [0, T]$
are related to the eigenvalue distributions
of random matrices in Gaussian ensembles
and in some two-matrix models.
Extreme values of paths in $[0, T]$ are
studied for these noncolliding diffusion processes
and determinantal and pfaffian representations
are given for the distribution functions.
The entries of the determinants and pfaffians
are expressed using special functions.
\footnote{2000 Mathematics Subject Classification(s):
Primary 15B52, 60J60, 62G32; Secondary 82C22, 17B10.}
\end{abstract}


\SSC{Introduction}

For $N \geq 2$, consider the following region in $\R^N$,
\begin{equation}
\W_N^{\rm A}=\Big\{
\x=(x_1, x_2, \dots, x_N) \in \R^N :
x_1 < x_2 < \cdots < x_N \Big\},
\label{eqn:WeylA}
\end{equation}
which is called the Weyl chamber of type ${\rm A}_{N-1}$
in the representation theory \cite{FH91}.
By the Karlin-McGregor formula \cite{KM59}, 
the density at $\y \in \W_N^{\rm A}$ 
of an $N$-dimensional Brownian motion
at time $t>0$, which starts from $\x \in \W_N^{\rm A}$
at time 0, and is 
restricted on the event that it stays in $\W_N^{\rm A}$
during a time interval $[0, t]$, is given by
$$
f^{\rm A}_{N}(t, \y|\x)= \det_{1 \leq i, j \leq N}
\Big[ p(t, y_i|x_j) \Big],
$$
where $p(t, y|x)$ is the heat kernel, 
$p(t,y|x)=e^{-(x-y)^2/(2t)}/\sqrt{2 \pi t}$.
It can be regarded as the transition probability density
of the absorbing Brownian motion in $\W_N^{\rm A}$
from $\x \in \W_N^{\rm A}$ to $\y \in \W_N^{\rm A}$
with duration $t >0$.

Set $0 < T < \infty$, $\a=(a_1, \dots, a_N), 
\b=(b_1, \dots, b_N) \in \W_N^{\rm A}$.
Then we consider the $N$-particle system of one-dimensional
Brownian motions starting from a configuration
$\a$ at time 0 and arriving at the configuration $\b$
at time $T$, conditioned never to collide with each other
in $[0, T]$.
The probability density at the configuration $\x \in \W_N^{\rm A}$
at time $t \in [0, T]$ is given by
\begin{equation}
p_{N, T}^{\rm A}(t, \x; \a, \b)
=\frac{f^{\rm A}_N(T-t, \b|\x) f^{\rm A}_N(t, \x|\a)}
{f^{\rm A}_N(T, \b|\a)}.
\label{eqn:pA1}
\end{equation}
We will call the above mentioned process 
the {\it noncolliding Brownian motion
from $\a$ to $\b$ with duration $T$}
and write it as 
${\X^{\a, \b}}(t)=(X^{\a, \b}_1(t), 
\dots, X^{\a, \b}_N(t)), t \in [0, T]$.

We can show that the limit
of (\ref{eqn:pA1}) in $\a \to \0 \equiv (0,0, \dots, 0)$,
$\b \to \0$, is given by
\begin{eqnarray}
\label{eqn:pA2}
\lim_{|\a|, |\b| \to 0}
p_{N, T}^{\rm A}(t, \x; \a, \b)
&=& p_{N, T}^{\rm A}(t, \x; \0, \0) \\
&=& q_{N}^{\rm GUE}(\x; \sigma_T(t)^2),
\quad t \in [0, T], \quad
\x \in \W_N^{\rm A}
\nonumber
\end{eqnarray}
with $\sigma_T(t)=\sqrt{t(1-t/T)}$ 
(Proposition 13 in \cite{KT04}).
Here $q_N^{\rm GUE}(\x; \sigma^2)$
denotes the probability density of eigenvalues
$\x \in \W_N^{\rm A}$ of random matrices in the
Gaussian unitary ensemble (GUE) with variance $\sigma^2$,
$$
q_{N}^{\rm GUE}(\x; \sigma^2)
= \frac{\sigma^{-N^2}}{(2 \pi)^{N/2}
\prod_{k=1}^{N} \Gamma(k)}
e^{-|\x|^2/(2 \sigma^2)}
\prod_{1 \leq i < j \leq N}
(x_j-x_i)^2,
$$
where $|\x|^2=\sum_{j=1}^{N} x_j^2$
and $\Gamma(k)$ is the Gamma function.
Note that, when $k \in \N \equiv \{1,2, \dots\}$, 
$\Gamma(k)=(k-1)!$.
We regard (\ref{eqn:pA2}) as 
the probability density at a configuration $\x \in \W_N^{\rm A}$
at time $t \in [0, T]$ of the
{\it noncolliding Brownian bridges with duration $T$}
denoted by 
$\X^{\0, \0}(t)=(X^{\0, \0}_1(t), 
\dots, X^{\0, \0}_N(t)), t \in [0, T]$.

When we make the configuration $\b$ at time $T$
be arbitrary in $\W_N^{\rm A}$, we have another
temporally inhomogeneous system of noncolliding Brownian motion,
which is denoted by 
$\X^{\a, \R}(t)=(X^{\a, \R}_1(t), 
\dots, X^{\a, \R}_N(t)), t \in [0, T]$.
The probability density at $\x \in \W_N^{\rm A}$
at time $t \in [0, T]$ is given by
$$
p_{N,T}^{\rm A}(t, \x; \a, \R)
=\frac{\cN^{\rm A}_{N} (T-t, \x) f^{\rm A}_{N} (t, \x|\a)}
{\cN^{\rm A}_N(T, \a)},
\quad t \in [0, T], 
$$
where
$$
\cN^{\rm A}_N(s, \x)=\int_{\W_N^{\rm A}} d \y \,
f^{\rm A}_N(s, \y|\x), \quad
s > 0, \quad \x \in \W_N^{\rm A}, 
$$
is the probability that the absorbing Brownian motion
in $\W_N^{\rm A}$ starting from $\x \in \W_N^{\rm A}$
is not yet absorbed at any boundary of the region
$\W_N^{\rm A}$ and is surviving inside of it at time $s >0$.
For an even integer $n$ and an antisymmetric $n\times n$ matrix 
$A = (a_{ij})$ we put
$$
\Pf_{1\le i, j \le n} \Big[a_{ij} \Big]
= \frac{1}{(n/2)!} 
\sum_{\sigma:\sigma(2k-1)<\sigma(2k), 1 \leq k \leq n/2} 
{\rm sgn} (\sigma)
a_{\sigma(1)\sigma(2)}a_{\sigma(3)\sigma(4)}
\cdots a_{\sigma(n-1)\sigma(n)},
$$
where the summation is extended over all permutations $\sigma$
of $(1,2,\dots,n)$ with restriction
$\sigma(2k-1)<\sigma(2k)$, $k=1,2,\dots,n/2$.
This expression is known as a pfaffian
(see, for example, \cite{Stem90}).
By using the de Bruijn identity \cite{deBr55}, 
which will be given as Lemma \ref{thm:lemma_deBruijn}
in Section 3, 
we have the formula \cite{KT02,KT03a}
$$
\cN^{\rm A}_N(s, \x) = \Pf_{1\le i, j\le \widehat{N}}
      \left[ F^{\rm A}_{ij} \left( \frac{\x}{\sqrt{2s}} \right) 
      \right],
$$
where
$$
\widehat{N} = \left\{
   \begin{array}{ll}
      N, & \quad \mbox{if} \ N=\mbox{even}, \\
      N+1, & \quad \mbox{if} \ N=\mbox{odd}, \\
   \end{array}
\right.
$$
$\x/\sqrt{2s}=(x_1/\sqrt{2s}, \dots, x_N/\sqrt{2s})$,
and
$$
F^{\rm A}_{ij}(\y) =
\left\{
   \begin{array}{ll}
\Psi(y_j-y_i),
      & \quad \mbox{if} \ 1\le i, j \le N,
\\
1, & \quad \mbox{if} \ 1\le i \le N, j=N+1,
\\
-1, & \quad \mbox{if} \ i=N+1, 1\le j \le N,
\\
0, & \quad \mbox{if} \ i=j=N+1,
\\
   \end{array}
\right.
$$
with
\begin{equation}
\Psi(u)= \frac{2}{\sqrt{\pi}} \int_0^u e^{-v^2}dv.
\label{eqn:PsiA}
\end{equation}

In \cite{KT02,KT03a}, the process
$\X^{\0, \R}(t) \equiv \lim_{|\a| \to 0} \X^{\a, \R}(t)
=(X^{\0, \R}_1(t),
\dots, X^{\0, \R}_N(t)),
t \in [0, T]$ is constructed, where 
the probability density at $\x \in \W_N^{\rm A}$
at time $t \in [0, T]$ is given by
\begin{equation}
p_{N, T}^{\rm A}(t, \x; \0, \R)
=\frac{T^{N(N-1)/4} t^{-N^2/2}}
{2^{N/2} \prod_{k=1}^{N} \Gamma(k/2)}
e^{-|\x|^2/(2t)} 
\prod_{1 \leq i < j \leq N} (x_j-x_i)
\cN(T-t, \x).
\label{eqn:pA4}
\end{equation}
It has been shown that (\ref{eqn:pA4}) exhibits 
a transition from the eigenvalue distribution
of GUE to the eigenvalue distribution
of another ensemble of random matrices called
the Gaussian orthogonal ensemble (GOE)
as $t \to T$ \cite{KT03a,KT03b}.
In other words, there establishes an interesting
correspondence \cite{KT02} between the temporally inhomogeneous
system of noncolliding Brownian motion
$\X^{\0, \R}(t), t \in [0, T]$
and the two-matrix model of Pandey and Mehta \cite{PM83}
in the random matrix theory \cite{Mehta04}.

\begin{figure}[htbp]
  \begin{center}
   \includegraphics[width=1.00\linewidth]{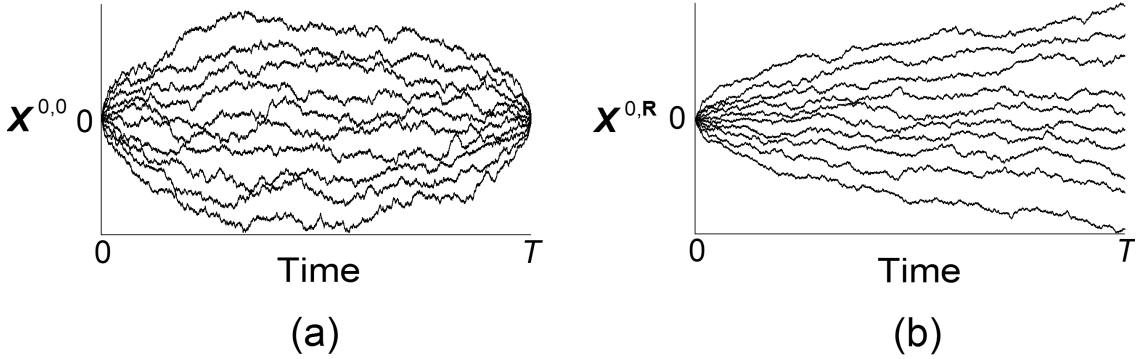}
  \end{center}
  \caption{Samples of paths for
  (a) $\X^{\0,0}(t)$
  and (b) $\X^{\0, \R}(t), t \in [0, T]$,
  generated by simulating the corresponding
  eigenvalue processes of random-matrix models.}
  \label{fig:fig_1}
\end{figure}

The equivalence of the distribution of $\X^{\0, \0}(t)$
and the eigenvalue distribution 
of random matrices in GUE with variance $\sigma_T(t)^2$
for $t \in [0, T]$ shown by Eq.(\ref{eqn:pA2}),
and the correspondence between $\X^{\0, \R}(t),
t \in [0, T]$ and the two-matrix model
of Pandey and Mehta are very useful
to perform computer simulations of the conditional
diffusion processes $\X^{\0, \0}(t)$ and
$\X^{\0, \R}(t), t \in [0, T]$ \cite{KIK08b}.
Figure \ref{fig:fig_1} shows samples of paths
generated by computer simulations
for these two processes.

In the present paper, we study the distributions
of the extreme values defined by
\begin{eqnarray}
\label{eqn:LRA} 
L^{\sharp}(N, T) &=& \min_{t \in [0, T]}
\X^{\sharp}(t) = \min_{t \in [0, T]} X^{\sharp}_1(t),
\\
R^{\sharp}(N, T) &=& \max_{t \in [0, T]}
\X^{\sharp}(t) = \max_{t \in [0, T]} X^{\sharp}_N(t),
\nonumber
\end{eqnarray}
for $\sharp=(\a, \b) \in \W_N^{\rm A} \times \W_N^{\rm A}$,
$\sharp=(\0, \0)$, and $(\0, \R)$.

In this paper, we also consider the
noncolliding processes ``with a wall at the origin".
Instead of the Weyl chamber of type ${\rm A}_{N-1}$
given by (\ref{eqn:WeylA}), we consider the
Weyl chamber of type ${\rm C}_N$,
$$
\W_N^{\rm C}=\Big\{
\x=(x_1, x_2, \dots, x_N) \in \R^N :
0< x_1 < x_2 < \cdots < x_N \Big\}.
$$
The density at $\y \in \W_N^{\rm C}$ 
of an $N$-dimensional Brownian motion
at time $t>0$, which starts from $\x \in \W_N^{\rm C}$
at time 0, and is
restricted on the event that it stays in $\W_N^{\rm C}$
during a time interval $[0, t]$, is given by
the Karlin-McGregor formula as
\begin{equation}
f^{\rm C}_{N}(t, \y|\x)= \det_{1 \leq i, j \leq N}
\Big[ p_{\rm abs}^{(0, \infty)}(t, y_i|x_j) \Big],
\label{eqn:fNC}
\end{equation}
where $p_{\rm abs}^{(0, \infty)}(t, y|x)$ is the transition
probability density of the one-dimensional
absorbing Brownian motion with an absorbing wall at the
origin. By the reflection principle of
Brownian motion, it is given as
\begin{eqnarray}
p_{\rm abs}^{(0, \infty)}(t, y|x) &=&
p(t, y|x)-p(t, y|-x) 
\nonumber\\
&=& \frac{1}{\sqrt{ 2 \pi t}}
\Big\{ e^{-(x-y)^2/(2t)} -e^{-(x+y)^2/(2t)} \Big\}
\nonumber
\end{eqnarray}
Eq.(\ref{eqn:fNC}) gives
the transition probability density
of the absorbing Brownian motion in $\W_N^{\rm C}$
from $\x$ to $\y$ with duration $t >0$.

For $0 < T < \infty$, $\a, \b \in \W_N^{\rm C}$,
we consider the $N$-particle system of one-dimensional
Brownian motions starting from a configuration
$\a$ at time 0 and arriving at a configuration $\b$
at time $T$, which is conditioned so that 
there is no collision between any pair of particles
and that all particles stay positive in a
time interval $[0, T]$.
The probability density at $\x \in \W_N^{\rm C}$
at time $t \in [0, T]$ is given by
\begin{equation}
p_{N, T}^{\rm C}(t, \x; \a, \b)
=\frac{f^{\rm C}_N(T-t, \b|\x) f^{\rm C}_N(t, \x|\a)}
{f^{\rm C}_N(T, \b|\a)}.
\label{eqn:pC1}
\end{equation}
The Brownian motion conditioned to stay positive
for a finite time interval $[0, T]$
is called the Brownian meander with duration $T$
\cite{RY98}.
So we call the $N$-particle system, whose 
probability density is given by (\ref{eqn:pC1}), 
the {\it noncolliding Brownian meander
from $\a$ to $\b$ with duration $T$} \cite{KT07a}
and write it as 
$\Y^{\a, \b}(t)=(Y^{\a, \b}_1(t), 
\dots, Y^{\a, \b}_N(t)), t \in [0, T]$.

The three-dimensional Bessel bridge with
duration $T$ is the conditional
Brownian motion such that 
it starts from the origin, stays positive in $t \in (0, T)$,
and first returns to the origin at time $T$.
Then the $N$-particle system obtained by the limit
$\Y^{\0, \0}(t) \equiv \lim_{|\a|, |\b| \to 0} \Y^{\a, \b}(t)
=(Y^{\0, \0}_1(t), 
\dots, Y^{\0, \0}_N(t)),
t \in [0, T]$ can be called
the {\it noncolliding three-dimensional Bessel bridge}.
It was called the system of
nonintersecting Brownian excursions in \cite{TW07}. 
The probability density at $\x \in \W_N^{\rm C}$
at time $t \in [0, T]$ is given by
$$
p_{N, T}^{\rm C}(t, \x; \0, \0)
= q_{N}^{\rm C}(\x; \sigma_T(t)^2)
$$
with $\sigma_T(t)=\sqrt{t(1-t/T)}$ 
(Proposition 14 in \cite{KT04}).
Here $q_N^{\rm C}(\x; \sigma^2)$
denotes the probability density of 
positive eigenvalues
$\x \in \W_N^{\rm C}$ of random matrices in the ensemble called
the class {\rm C} with variance $\sigma^2$
\cite{AZ97},
$$
q_{N}^{\rm C}(\x; \sigma^2)
= \frac{\sigma^{-N(2N+1)}}{(\pi/2)^{N/2}
\prod_{\ell=1}^{N} \Gamma(2 \ell)}
e^{-|\x|^2/(2 \sigma^2)}
\prod_{1 \leq i < j \leq N}
(x_j^2-x_i^2)^2 \prod_{k=1}^{N} x_k^2.
$$

When we make the configuration $\b$ at time $T$
be arbitrary in $\W_N^{\rm C}$, we have another
noncolliding Brownian meander,
which is denoted by \\
$\Y^{\a, \R_+}(t)=(Y^{\a, \R_+}_1(t),
\dots, Y^{\a, \R_+}_N(t))$, $t \in [0, T]$
with $\R_+=\{x \in \R: x >0\}$.
The probability density at $\x \in \W_N^{\rm C}$
at time $t \in [0, T]$ is given by
$$
p_{N,T}^{\rm C}(t, \x; \a, \R_+)
=\frac{\cN^{\rm C}_N(T-t, \x) f^{\rm C}_N(t, \x|\a)}
{\cN^{\rm C}_N(T, \a)},
$$
where
$$
\cN^{\rm C}_N(s, \x)=\int_{\W_N^{\rm C}} d \y \,
f^{\rm C}_N(s, \y|\x), \quad
s > 0, \quad \x \in \W_N^{\rm C}, 
$$
is the probability that the absorbing Brownian motion
in $\W_N^{\rm C}$ starting from $\x \in \W_N^{\rm C}$
is not yet absorbed at any boundary of
$\W_N^{\rm C}$ and is surviving inside of it at time $s >0$.
We can prove the formula \cite{KT02,KT03a}
$$
\cN^{\rm C}_N(s,\x) =
   \Pf_{1\le i, j\le \widehat{N}}
 \left[ F^{\rm C}_{ij}\left( \frac{\x}{\sqrt{2s}} \right) \right],
$$
where
$\x/\sqrt{2s}=(x_1/\sqrt{2s}, \dots, x_N/\sqrt{2s})$,
$$
F^{\rm C}_{ij}(\y) =
\left\{
   \begin{array}{ll}
   \Psi(y_i, y_j),
      & \quad \mbox{if} \ 1\le i, j \le N,
\\
\Psi(y_i),
 & \quad \mbox{if} \ 1\le i \le N, j=N+1,
\\
-\Psi(y_j),
& \quad \mbox{if} \ i=N+1, 1\le j \le N,
\\
0, & \quad \mbox{if} \ i=j=N+1,
\\
   \end{array}
\right.
$$
with (\ref{eqn:PsiA}) and
\begin{eqnarray}
\Psi(u_1,u_2)
&=& \frac{2}{\pi}
\left[ \int_0^{u_1} dv_1 \int_{u_1 -u_2}^{u_2 -u_1}dv_2 \,
e^{-v_1^2 - (v_1-v_2)^2}
\right.
\nonumber
\\
&&\quad - \left. 
\int_{u_1}^{u_2} dv_1 \int_{u_2 -u_1}^{u_1 +u_2}dv_2 \, 
e^{-v_1^2 - (v_1-v_2)^2}
\right].
\nonumber
\end{eqnarray}

The process
$\Y^{\0, \R_+}(t) \equiv \lim_{|\a| \to 0} \Y^{\a, \R_+}(t)
=(Y^{\0, \R_+}_1(t), 
\dots, Y^{\0, \R_+}_N(t)),
t \in [0, T]$ is studied in \cite{KT03a,KTNK03}, where 
the probability density at $\x \in \W_N^{\rm C}$
at time $t \in [0, T]$ is given by
\begin{eqnarray}
\label{eqn:pC4}
&& p_{N, T}^{\rm C}(t, \x; \0, \R_+) \\
&& \quad
=\frac{T^{N^2/2} t^{-N(2N+1)/2}}
{\prod_{\ell=1}^{N} \Gamma(\ell)}
e^{-|\x|^2/(2t)} 
\prod_{1 \leq i < j \leq N} (x_j^2-x_i^2) \prod_{k=1}^{N} x_k
\, 
\cN^{\rm C}(T-t, \x).
\nonumber
\end{eqnarray}
It has been shown that (\ref{eqn:pC4}) exhibits 
a transition from the eigenvalue distribution
of random matrices in the ensemble called the class {\rm C}
to the eigenvalue distribution
of the class {\rm CI}
as $t \to T$ \cite{KTNK03,KT04}.

\begin{figure}[htbp]
  \begin{center}
   \includegraphics[width=1.00\linewidth]{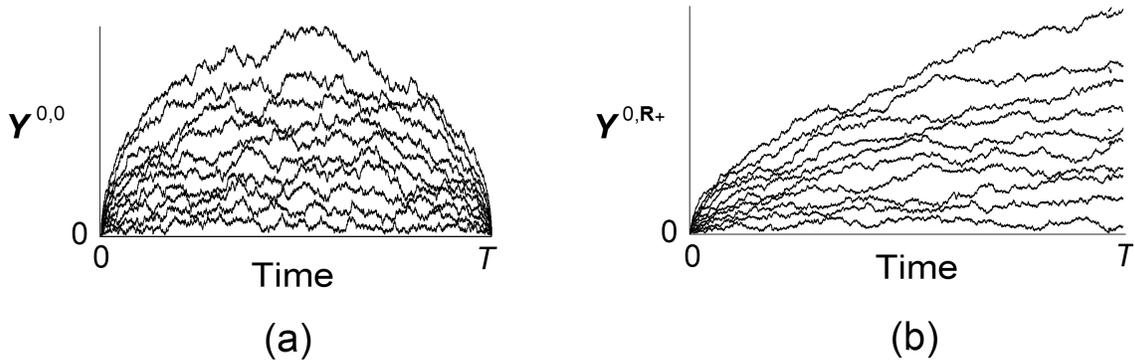}
  \end{center}
  \caption{Samples of paths for 
  (a) $\Y^{\0, \0}(t)$ and (b) $\Y^{\0, \R_+}(t), t \in [0, T]$.}
  \label{fig:fig_2}
\end{figure}

We can perform computer simulations of a matrix-valued
Brownian bridge such that its distribution is
the same as the distribution of random matrices
in the class {\rm C} with variance
$\sigma_T(t)=\sqrt{t(1-t/T)}, t \in [0, T]$,
and a matrix-valued Brownian motion,
whose distribution at time $t \in [0, T]$
changes continuously from the class {\rm C}
distribution to the class {\rm CI} distribution
of random matrices as $t \to T$.
By numerically calculating eigenvalues
of these two matrix-valued diffusion processes
in $t \in [0, T]$, we can draw samples of paths
of noncolliding diffusion particles
for $\Y^{\0, \0}(t)$ and $\Y^{\0, \R_+}(t),
t \in [0, T]$, as shown in Figure \ref{fig:fig_2}
\cite{KIK08b}.

For $\Y^{\sharp}$, 
$\sharp=(\a, \b) \in \W_N^{\rm C} \times \W_N^{\rm C}$, 
$(\0, \0)$, and $(\0, \R_+)$,
we will report the distributions
of the maximum value
\begin{equation}
H^{\sharp}(N, T)= \max_{t \in [0, T]} \Y^{\sharp}(t)
=\max_{t \in [0, T]} Y^{\sharp}_N(t).
\label{eqn:HC}
\end{equation}

\SSC{Results}

For $k \in \N_0 \equiv \{0,1,2, \dots\}$, let
$H_k(x)$ be the $k$-th Hermite polynomial
$$
H_k(x)=k ! \sum_{m=0}^{[k/2]} 
\frac{(-1)^m (2x)^{k-2m}}{m ! (k-2m)!},
$$
where $[a]$ denotes the largest integer that is not
greater than $a$.
We define the function of $u, v \in \R$ with
an index $k \in \N_0$,
\begin{equation}
\Theta_k(u, v)
=\sum_{n \in \Z} H_k(un+v) e^{-(un+v)^2}.
\label{eqn:Thetak}
\end{equation}
If we consider the following version of Jacobi's theta function
$$
\vartheta(x,y)=\sum_{n \in \Z}
q^{n^2} z^{2n}
=\sum_{n \in \Z}
e^{2 \pi \sqrt{-1} xn +\pi \sqrt{-1} y n^2},
$$
${\rm Im} \, y >0$, where we have set
$z=e^{\pi \sqrt{-1} x}$ and $q=e^{\pi \sqrt{-1} y}$,
the reciprocal relation
$$
\vartheta(x,y)=\vartheta
\left( \frac{x}{y}, - \frac{1}{y} \right)
e^{-\pi \sqrt{-1} x^2/y} 
\left( \frac{\sqrt{-1}}{y} \right)^{1/2}
$$
holds (see Sec.10.12 in \cite{AAR99},
Sec.A.3.1 in \cite{Ful07}).
Combining this functional equation with the fact
$$
e^{2ux-u^2}=\sum_{k=0}^{\infty} H_k(x)
\frac{u^k}{k!},
$$
the following identities are derived
for $k \in \N_0, \eta, \xi \in \R$
(Lemma 3 in \cite{KIK08b});
\begin{equation}
\sum_{n \in \Z} n^k
\exp \left(-\frac{\pi}{\eta^2} n^2
+\frac{2 \pi \sqrt{-1} \xi}{\eta^2} n \right)
= \frac{(-1)^{k/2} \eta^{k+1}}
{2^k \pi^{k/2}}
\Theta_k \left( \sqrt{\pi} \eta,
\frac{\sqrt{\pi} \xi}{\eta} \right).
\label{eqn:Iden1}
\end{equation}

The following two theorems are our main results.

\begin{thm}
\label{thm:Thm2_1}
For $\ell, r > 0$,
\begin{eqnarray}
\label{eqn:thm2_1_1}
&& \P \Big( -\ell < L^{\0, \0}(N,T),
R^{\0, \0}(N,T) < r \Big)
\\
&& = \frac{(-1)^{N}}{2^{N(N-1)/2} \prod_{k=1}^{N} \Gamma(k)}
\nonumber\\
&& \times
\det_{1 \leq i, j \leq N}
\left[ \Theta_{i+j-2}
\left( \frac{2(\ell+r)}{\sqrt{2T}},
\frac{2 \ell}{\sqrt{2T}} \right)
+(-1)^{j}
\Theta_{i+j-2}
\left( \frac{2(\ell+r)}{\sqrt{2T}}, 0 \right)
\right],
\nonumber
\end{eqnarray}
and
\begin{eqnarray}
\label{eqn:thm2_1_2}
&& \P \Big( -\ell < L^{\0, \R}(N,T),
R^{\0, \R}(N,T) < r \Big)
\\
&& = \frac{(-1)^{N(N+1)/2}}
{2^{N(N-1)/4} \prod_{k=1}^{N} \Gamma(k/2)} 
\Pf_{1 \leq i, j \leq \widehat{N}}
\left[ G_{ij}^{\rm A} \left(-\frac{\ell}{\sqrt{T}},
\frac{r}{\sqrt{T}}\right) \right],
\nonumber
\end{eqnarray}
where
\begin{equation}
G_{ij}^{\rm A} \left(-\frac{\ell}{\sqrt{T}},
\frac{r}{\sqrt{T}}\right)
= \left\{ \begin{array}{rl}
I_{ij}^{\rm A}, \quad & \mbox{if} \quad 1 \leq i, j \leq N,
\cr
I_{i}^{\rm A}, \quad & \mbox{if} \quad 
1 \leq i \leq N, j=N+1, 
\cr
-I_{j}^{\rm A}, \quad & \mbox{if} \quad
i=N+1, 1 \leq j \leq N,
\cr
0, \quad & \mbox{if} \quad
i=j=N+1,
\end{array} \right.
\nonumber
\end{equation}
with
\begin{eqnarray}
&& z_i^{\rm A}(x)= 
\Theta_{i-1} \left( \frac{2(\ell+r)}{\sqrt{2T}},
\frac{2\ell}{\sqrt{2T}}+x \right)
+(-1)^{i} 
\Theta_{i-1} \left( \frac{2(\ell+r)}{\sqrt{2T}}, x \right)
\nonumber\\
&& I_i^{\rm A} =
\int_{-\ell/\sqrt{2T}}^{r/\sqrt{2T}} 
z_i^{\rm A}(x) dx,
\nonumber\\
&& I_{ij}^{\rm A} =
\int_{-\ell/\sqrt{2T} < x_1 < x_2 < r/\sqrt{2T}} 
\det \left[ \begin{array}{ll}
z_i^{\rm A}(x_1) & z_i^{\rm A}(x_2) \cr
z_j^{\rm A}(x_1) & z_j^{\rm A}(x_2) 
\end{array} \right] dx_1 dx_2.
\nonumber
\end{eqnarray}
\end{thm}

\begin{thm}
\label{thm:Thm2_2}
For $h > 0$,
\begin{equation}
\label{eqn:thm2_2_1}
\P \Big(H^{\0, \0}(N,T) < h  \Big)
= \frac{(-1)^N}{2^{N^2} 
\prod_{k=1}^{N} \Gamma(2k)}
\det_{1 \leq i, j \leq N} \left[
\Theta_{2(i+j-1)}
\left( \frac{2h}{\sqrt{2T}}, 0 \right) \right],
\end{equation}
and
\begin{equation}
\label{eqn:thm2_2_2}
\P \Big( H^{\0, \R_+}(N,T) < h \Big)
= \frac{1}{2^{N(N-1)/2} \prod_{k=1}^{N} \Gamma(k)}
\Pf_{1\le i, j\le \widehat{N}}
\left[ G_{ij}^{\rm C} \left(\frac{h}{\sqrt{T}} \right) \right],
\end{equation}
where
\begin{equation}
G_{ij}^{\rm C} \left( \frac{h}{\sqrt{T}} \right)
= \left\{ \begin{array}{rl}
I_{ij}^{\rm C}, \quad & \mbox{if} \quad 1 \leq i, j \leq N,
\cr
I_{i}^{\rm C}, \quad & \mbox{if} \quad 
1 \leq i \leq N, j=N+1, 
\cr
-I_{j}^{\rm C}, \quad & \mbox{if} \quad
i=N+1, 1 \leq j \leq N,
\cr
0, \quad & \mbox{if} \quad
i=j=N+1,
\end{array} \right.
\nonumber
\end{equation}
with
\begin{eqnarray}
&& z_i^{\rm C}(x)= 
\Theta_{2i-1} \left( \frac{2h}{\sqrt{2T}}, x \right)
\nonumber\\
&& I_i^{\rm C} =
\int_{0}^{h/\sqrt{2T}} 
z_i^{\rm C}(x) dx,
\nonumber\\
&& I_{ij}^{\rm C} =
\int_{0 < x_1 < x_2 < h/\sqrt{2T}} 
\det \left[ \begin{array}{ll}
z_i^{\rm C}(x_1) & z_i^{\rm C}(x_2) \cr
z_j^{\rm C}(x_1) & z_j^{\rm C}(x_2) 
\end{array} \right] dx_1 dx_2.
\nonumber
\end{eqnarray}
\end{thm}

\noindent{\bf Remarks.} 
\begin{description}
\item{(1)} \quad
By the scaling property of Brownian motion,
we have the equalities in distribution
$$
L^{\sharp}(N,T) \=
\sqrt{T} L^{\sharp}(N,1), \quad
 R^{\sharp}(N,T) \=
\sqrt{T} R^{\sharp}(N,1),
$$
for $\sharp=(\0,\0), (\0, \R)$,
and
$$
H^{\sharp}(N,T)=\sqrt{T} H^{\sharp}(N,1)
$$
for $\sharp=(\0, \0), (\0, \R_+)$, $\forall T>0$.
As a matter of fact,
the probability distributions (\ref{eqn:thm2_1_1}) and
(\ref{eqn:thm2_1_2}) in Theorem \ref{thm:Thm2_1}
are functions of $\ell/\sqrt{T}$ and $r/\sqrt{T}$,
and (\ref{eqn:thm2_2_1}) and (\ref{eqn:thm2_2_2})
are of $h/\sqrt{T}$, respectively.

\item{(2)} \quad
For $N=1$, (\ref{eqn:thm2_2_1}) gives
\begin{eqnarray}
\P\Big( H^{\0, \0}(1, T) < h \Big)
&=& - \frac{1}{2} \sum_{n \in \Z} 
H_2 \left( \frac{2h}{\sqrt{2T}} n \right)
e^{-2h^2n^2/T}
\nonumber\\
&=& \sum_{n \in Z} \left(1-\frac{4h^2n^2}{T} \right)
e^{-2h^2 n^2/T},
\nonumber
\end{eqnarray}
since $H_2(x)=4x^2-2$.
This is a classical result for the height distribution
of the three-dimensional Bessel bridge.
It should be remarked that this result for
distribution is equivalent with the fact on moments
$$
\E \Big[
\Big( H^{\0, \0}(1,T) \Big)^m \Big]
=2 \left( \frac{\pi T}{2} \right)^{m/2}
\xi(m), \quad m \in \C
$$
with
$$
\xi(m)=\frac{1}{2} m(m-1) \pi^{-m/2} 
\Gamma \left(\frac{m}{2} \right) \zeta(m),
$$
where $\zeta(m)$ is the Riemann zeta function
$$
\zeta(m)=\sum_{n \in \N} \frac{1}{n^m},
\quad {\rm Re} \, m > 1.
$$
See \cite{BPY01} for interesting relations
of the probability laws of one-dimensional conditional
Brownian motions
to the Jacobi theta functions and the
Riemann zeta function.

\item{(3)} \quad
In \cite{KIK08a}, $N=2$ case of the formula (\ref{eqn:thm2_2_1})
of Theorem \ref{thm:Thm2_2}
and moments calculated from it have been intensively studied,
which are related to the double Dirichlet series.

\item{(4)} \quad
The formula (\ref{eqn:thm2_1_1}) of Theorem \ref{thm:Thm2_1}
gives the joint distribution of two extreme values 
$L^{\0,\0}(N,T)$ and $R^{\0, \0}(N,T)$
of the process $\X^{\0, \0}(t), t \in [0, T]$.
The distributions of single variable
$R^{\0, \0}(N,T)$ was studied in \cite{Feierl09}
and \cite{SMCRF08}.
Feierl \cite{Feierl09} also reported the distribution of
the width (or ``range") defined by
\begin{eqnarray}
W^{\0, \0}(N, T)
&=& R^{\0, \0}(N,T)-L^{\0, \0}(N,T)
\nonumber\\
&=& \max_{t \in [0, T]} \Big(
X_N^{\0, \0}(t)-X^{\0, \0}_1(t) \Big).
\nonumber
\end{eqnarray}

\item{(5)} \quad
A part of results recently reported 
by Borodin {\it et al.} \cite{BFPSW09}
may be stated as follows in the present 
notations;
if $N$ is even, for each fixed $T \geq 0$,
$$
\max_{t \in [0, T/2]}
\frac{1}{\sqrt{2(1-t/T)}} X^{\0,\0}_N(t) \= Y^{\0, \0}_{N/2}(T/2),
$$
where $Y^{\0, \0}_{N/2}(T/2)$ in the RHS denotes the
position of the top particle
at time $T/2$ of the $N/2$-particle
system $\Y^{\0, \0}(t)=(Y^{\0, \0}_1(t),
\dots, Y^{\0, \0}_{N/2}(t)), t \in [0, T]$.
The distribution of the top path 
$Y^{\0, \0}_{N/2}(t), t \in [0, T]$
was studied by Tracy and Widom \cite{TW07}.

\item{(6)} \quad
The distribution of the ``maximum height"
$H^{\0, \0}(N,T)$ of the process
$\Y^{\0,\0}(t), t \in [0, T]$
has been studied by Feierl \cite{Feierl07}
and by Schehr {\it et al.} \cite{SMCRF08}
in different context and by different methods.
See \cite{KIK08b} for comparison.
\end{description}

\SSC{Proof of Theorem \ref{thm:Thm2_1}}

Let ${\bf 1}\{\omega\}$ be the indicator function
of a condition $\omega$.
Assume that $- \ell < a_1, - \ell < b_1,
a_N < r, b_N < r$.
By definition of the process
and (\ref{eqn:LRA}),
for $(\a, \b) \in \W_N^{\rm A} \times \W_N^{\rm A}$
\begin{eqnarray}
\label{eqn:PEE}
&& \P \Big( -\ell < L^{\a, \b}(N, T), 
R^{\a, \b}(N, T) < r \Big)
\\
&=& \E \Big(
{\bf 1}\{ - \ell < X^{\a, \b}_1(t), 
X^{\a, \b}_N(t) < r,
\, \forall t \in [0, T]\} \Big)
\nonumber\\
&=& \E \Big(
{\bf 1} \{ -\ell < X^{\a, \b}_i(t) < r,
\, 1 \leq \forall i \leq N, \,
\forall t \in [0, T] \} \Big).
\nonumber
\end{eqnarray}

For $\ell, r >0$, 
we introduce the following restricted region in 
$\W_N^{\rm A}$,
$$
\W_N^{\rm A}(-\ell, r)
=\Big\{ \x=(x_1, \dots, x_N) \in \W_N^{\rm A}:
-\ell < x_1, x_N < r \Big\}.
$$
Then we consider the absorbing Brownian motion
in this region
starting from $\a \in \W_N^{\rm A}(-\ell, r)$.
Let the transition probability density of
a one-dimensional Brownian motion in an interval
$(-\ell, r)$, where two absorbing walls are put
at $x=-\ell$ and $x=r$, be
$p_{\rm abs}^{(-\ell, r)}(t, \cdot \, | \, \cdot), t \geq 0$.
Then, by the Karlin-McGregor formula, the probability density
that the absorbing Brownian motion in $\W_N^{\rm A}(-\ell, r)$
arrives at $\b \in \W_N^{\rm A}(-\ell, r)$
at time $T >0$, avoided from any absorption
at boundary of the region $\W_N^{\rm A}(-\ell, r)$,
is given by
$$
f_N^{A(-\ell, r)}(T, \b|\a)
=\det_{1 \leq i, j \leq N}
\Big[ p_{\rm abs}^{(-\ell, r)}(T, b_i|a_j) \Big].
$$
The probability density $p_{\rm abs}^{(-\ell, r)}(t, y|x)$
is obtained by solving the diffusion equation
$\partial u(t,y)/\partial t=(1/2) \partial^2 u(t,y)
/\partial y^2, -\ell \leq y \leq r$ with
the Dirichlet boundary condition
$u(\cdot \, , -\ell)=u(\cdot \, , r)=0$
and with the initial condition $u(0,y)=\delta_x(y),
-\ell < x, y < r$.
The solution is given by
$$
p_{\rm abs}^{(-\ell, r)}(t, y|x)
= \frac{2}{\ell+r} \sum_{n \in \N}
e^{-\pi^2 t n^2/\{2(\ell+r)^2\}}
\sin \left( \pi \frac{\ell+x}{\ell+r} n \right)
\sin \left( \pi \frac{\ell+y}{\ell+r} n \right),
$$
$t > 0, x,y \in (-\ell, r)$.

Then the following equality will be established,
\begin{equation}
\label{eqn:equal1}
\E \Big( {\bf 1} \{
-\ell < X_i^{\a, \b}(t) < r, \,
1 \leq \forall i \leq N, \,
\forall t \in [0, T] \} \Big)
=\frac{f_N^{\rm A(-\ell, r)}(T, \b|\a)}
{f_N^{\rm A}(T, \b|\a)}.
\end{equation}
And then we have the equalities
\begin{eqnarray}
\label{eqn:equal2}
&& \P \Big( -\ell < L^{\a, \R}(N,T),
R^{\a, \R}(N, T) < r \Big)
\\
&& \quad =
\E \Big( {\bf 1} \{
-\ell < X_i^{\a, \R}(t) < r, \,
1 \leq \forall i \leq N, \,
\forall t \in [0, T] \} \Big)
\nonumber\\
&& \quad
=\frac{ \displaystyle{
\int_{\W_N^{\rm A}(-\ell, r)} 
d \b \, f_N^{\rm A(-\ell, r)}(T, \b|\a)}}
{\cN^{\rm A}_N(T, \a)}.
\nonumber
\end{eqnarray}

Now we want to take the limit $\a, \b \to \0$ in
(\ref{eqn:equal1}) and $\a \to \0$ in (\ref{eqn:equal2})
in order to prove Theorem \ref{thm:Thm2_1}. 
The following asymptotics in $|\a| \to 0$
for the denominators of
(\ref{eqn:equal1}) and (\ref{eqn:equal2}) are
given in \cite{KT04,KT07b};
as $|\a| \to 0$, 
\begin{eqnarray}
\label{eqn:asym1}
&& f_N^{\rm A}(T, \b|\a)
= \frac{T^{-N^2/2}}{(2\pi)^{N/2} \prod_{k=1}^{N} \Gamma(k)}
e^{-|\b|^2/(2T)}
\\
&& \qquad \qquad \quad \times 
\prod_{1 \leq i < j \leq N} 
\Big\{(a_j-a_i)(b_j-b_i)\Big\}
\times \Big\{1+{\cal O}(|\a|) \Big\},
\nonumber\\
\label{eqn:asym2}
&& \cN_N^{\rm A}(T, \a)
= \frac{1}{\pi^{N/2}}
\prod_{k=1}^{N} \frac{\Gamma(k/2)}{\Gamma(k)}
T^{-N(N-1)/4} 
\\
&& \qquad \qquad \quad \times
\prod_{1 \leq i < j \leq N} (a_j-a_i)
\times \Big\{ 1+{\cal O}(|\a|) \Big\}.
\nonumber
\end{eqnarray}

By multilinearity of determinants, we have
\begin{eqnarray}
\label{eqn:fAlr0}
&& f_N^{{\rm A}(-\ell, r)}(T, \b|\a)
\\
&=& \det_{1 \leq i, j \leq N} \left[
\frac{2}{\ell+r} \sum_{n_i \in \N}
e^{-\pi^2 T n_i^2/\{2(\ell+r)^2\}}
\sin \left( \pi \frac{\ell+a_i}{\ell+r} n_i \right)
\sin \left( \pi \frac{\ell+b_j}{\ell+r} n_i \right) \right]
\nonumber\\
&=& \frac{1}{N!} \left(\frac{2}{\ell+r} \right)^{N}
\sum_{\n \in \N^N} 
e^{-\pi^2 T |\n|^2/\{2(\ell+r)^2\}}
\nonumber\\
&& \times \det_{1 \leq i, j \leq N}
\left[ \sin \left( \pi \frac{\ell+a_i}{\ell+r} n_j \right) \right]
\det_{1 \leq \alpha, \beta \leq N}
\left[ \sin \left( \pi \frac{\ell+b_{\alpha}}{\ell+r} 
n_{\beta} \right) \right],
\nonumber
\end{eqnarray}
where $|\n|^2=\sum_{i=1}^{N} n_i^2$.
We will use the series expansion
\begin{equation}
\sin \left( \pi \frac{\ell+a}{\ell+r} n \right)
=\sum_{p \in \N_0} c_p(n) (a n)^p
\label{eqn:sin1}
\end{equation}
with
\begin{equation}
c_p(n)=\frac{(-1)^{(p-1)/2}}{2 p!}
\left( \frac{\pi}{\ell+r} \right)^p
\left\{ \exp \left( \frac{\pi \sqrt{-1} \ell}{\ell+r} n \right)
+(-1)^{p+1}
\exp \left( -\frac{\pi \sqrt{-1} \ell}{\ell+r} n \right) \right\},
\label{eqn:cpn}
\end{equation}
$p \in \N_0$.
We will also use the property of the Vandermonde determinant
$$
\det_{1 \leq i, j \leq N} \Big[a_i^{j-1} \Big]
=\prod_{1 \leq i<j \leq N} (a_j-a_i).
$$
Then
\begin{eqnarray}
\label{eqn:fAlr0b}
&& f_{N}^{{\rm A}(-\ell, r)}(T, \b|\a)
\\
&=& \frac{1}{N!} \left( \frac{2}{\ell+r}\right)^N
\sum_{\n \in \N^{N}} 
e^{-\pi^2 T |\n|^2/\{2(\ell+r)^2\}}
\nonumber\\
&& \times \sum_{0 \leq p_1 < \cdots < p_N}
\det_{1 \leq i, j \leq N} \Big[ a_i^{p_j} \Big]
\det_{1 \leq \gamma, \delta \leq N}
\Big[ c_{p_{\gamma}}(n_{\delta}) n_{\delta}^{p_{\gamma}} \Big]
\det_{1 \leq \alpha, \beta \leq N}
\left[ \sin \left( \pi \frac{\ell+b_{\alpha}}{\ell+r}
n_{\beta} \right) \right]
\nonumber\\
&=& \frac{1}{N!} \left( \frac{2}{\ell+r}\right)^N
\sum_{\n \in \N^{N}} 
e^{-\pi^2 T |\n|^2/\{2(\ell+r)^2\}}
\det_{1 \leq \gamma, \delta \leq N}
\Big[ c_{\gamma-1}(n_{\delta}) n_{\delta}^{\gamma-1} \Big]
\nonumber\\
&& \times 
\det_{1 \leq \alpha, \beta \leq N}
\left[ \sin \left( \pi \frac{\ell+b_{\alpha}}{\ell+r}
n_{\beta} \right) \right]
\prod_{1 \leq i<j \leq N} (a_j-a_i)
\times \Big\{ 1+{\cal O}(|\a|) \Big\}
\nonumber\\
&=& \left( \frac{2}{\ell+r} \right)^N
\det_{1 \leq \alpha, \beta \leq N}
\left[ \sum_{n \in \N}
e^{-\pi^2 T n^2/\{2(\ell+r)^2\}}
c_{\alpha-1}(n) n^{\alpha-1}
\sin \left( \pi \frac{\ell+b_{\beta}}{\ell+r}
n \right) \right]
\nonumber\\
&& \qquad \times
\prod_{1 \leq i<j \leq N} (a_j-a_i)
\times \Big\{ 1+{\cal O}(|\a|) \Big\}.
\nonumber
\end{eqnarray}
By (\ref{eqn:cpn}), we have
\begin{eqnarray}
&& \sum_{n \in \N}
e^{-\pi^2 T n^2/\{2(\ell+r)^2\}}
c_{\alpha-1}(n) n^{\alpha-1}
\sin \left( \pi \frac{\ell+b_{\beta}}{\ell+r}
n \right)
\nonumber\\
&=& \frac{(-1)^{(\alpha-2)/2}}{2 (\alpha-1)!}
\left( \frac{\pi}{\ell+r} \right)^{\alpha-1}
\nonumber\\
&& \times \left[
\sum_{n \in \N} n^{\alpha-1}
\exp\left(-\frac{\pi^2 T}{2(\ell+r)^2} n^2
+\frac{\pi \sqrt{-1} \ell}{\ell+r} n \right)
\sin \left( \pi \frac{\ell+b_{\beta}}{\ell+r} n \right) \right.
\nonumber\\
&& \quad \left. +
\sum_{n \in \N} (-n)^{\alpha-1}
\exp\left(-\frac{\pi^2 T}{2(\ell+r)^2} (-n)^2
+\frac{\pi \sqrt{-1} \ell}{\ell+r} (-n) \right)
\sin \left( \pi \frac{\ell+b_{\beta}}{\ell+r} (-n) \right) \right]
\nonumber\\
&=& \frac{(-1)^{(\alpha-2)/2}}{2 (\alpha-1)!}
\left( \frac{\pi}{\ell+r} \right)^{\alpha-1}
\nonumber\\
&& \times
\sum_{n \in \Z} n^{\alpha-1}
\exp\left(-\frac{\pi^2 T}{2(\ell+r)^2} n^2
+\frac{\pi \sqrt{-1} \ell}{\ell+r} n \right)
\sin \left( \pi \frac{\ell+b_{\beta}}{\ell+r} n \right).
\nonumber\\
&=& \frac{(-1)^{(\alpha-3)/2}}{4 (\alpha-1)!}
\left( \frac{\pi}{\ell+r} \right)^{\alpha-1}
\left\{ \sum_{n \in \Z}
n^{\alpha-1} \exp \left(
-\frac{\pi^2 T}{2(\ell+r)^2} n^2
+\frac{\pi \sqrt{-1}(2\ell+b_{\beta})}{\ell+r} n \right)
\right.
\nonumber\\
&& \qquad \qquad  \left.
+(-1)^{\alpha}
\sum_{n \in \Z}
n^{\alpha-1} \exp \left(
-\frac{\pi^2 T}{2(\ell+r)^2} n^2
+\frac{\pi \sqrt{-1} b_{\beta}}{\ell+r} n \right)
\right\}.
\nonumber
\end{eqnarray}
By the identity (\ref{eqn:Iden1}),
the above is equal to 
\begin{eqnarray}
&& \frac{(-1)^{\alpha} (\ell+r)}
{2^{(\alpha+2)/2} \pi^{1/2} (\alpha-1)!} T^{-\alpha/2}
\nonumber\\
&& \quad \times
\left\{ \Theta_{\alpha-1}
\left( \frac{2(\ell+r)}{\sqrt{2T}},
\frac{2 \ell+b_{\beta}}{\sqrt{2T}} \right)
+(-1)^{\alpha}
 \Theta_{\alpha-1}
\left( \frac{2(\ell+r)}{\sqrt{2T}},
\frac{b_{\beta}}{\sqrt{2T}} \right) \right\}.
\nonumber
\end{eqnarray}
Then we have the following asymptotics in $|\a| \to 0$,
\begin{eqnarray}
\label{eqn:asym3}
&& f_N^{{\rm A}(-\ell, r)} (T, \b|\a)
\\
&=& \frac{(-1)^{N(N+1)/2}}
{2^{N(N+1)/4} \pi^{N/2}
\prod_{k=1}^{N} \Gamma(k)}
T^{-N(N+1)/4}
\nonumber\\
&&
\times \det_{1 \leq \alpha, \beta \leq N}
\left[\Theta_{\alpha-1}
\left( \frac{2(\ell+r)}{\sqrt{2T}},
\frac{2 \ell+b_{\beta}}{\sqrt{2T}} \right)
+(-1)^{\alpha}
 \Theta_{\alpha-1}
\left( \frac{2(\ell+r)}{\sqrt{2T}},
\frac{b_{\beta}}{\sqrt{2T}} \right) \right]
\nonumber\\
&& \qquad \times
\prod_{1 \leq i < j \leq N} (a_j-a_i)
\times \Big\{1+{\cal O}(|\a|) \Big\}.
\nonumber
\end{eqnarray}

From (\ref{eqn:fAlr0}) to (\ref{eqn:fAlr0b}),
we performed the series expansion (\ref{eqn:sin1})
in only one of the two determinants
in (\ref{eqn:fAlr0}).
If we perform the series expansions
in both of them, we have the following estimate,
\begin{eqnarray}
&& f_{N}^{{\rm A}(-\ell, r)}(T, \b|\a)
\nonumber\\
&=& \frac{1}{N!} \left( \frac{2}{\ell+r} \right)^N
\sum_{\n \in \N^N}
e^{-\pi^2 T |\n|^2/\{2(\ell+r)^2\}}
\nonumber\\
&& \times 
\sum_{0 \leq p_1 < \cdots < p_N}
\det_{1 \leq i, j \leq N}
\Big[ a_i^{p_j} \Big]
\det_{1 \leq \alpha, \beta \leq N}
\Big[ c_{p_{\alpha}}(n_{\beta})
n_{\beta}^{p_{\alpha}} \Big]
\nonumber\\
&& \times 
\sum_{0 \leq q_1 < \cdots < q_N}
\det_{1 \leq i', j' \leq N}
\Big[ b_{i'}^{q_{j'}} \Big]
\det_{1 \leq \alpha', \beta' \leq N}
\Big[ c_{q_{\alpha'}}(n_{\beta'})
n_{\beta'}^{q_{\alpha'}} \Big]
\nonumber\\
&=& \left(\frac{2}{\ell+r} \right)^N
\det_{1 \leq \alpha, \beta \leq N}
\left[ \sum_{n \in \N}
e^{-\pi^2 T n^2/\{2(\ell+r)^2\}}
c_{\alpha-1}(n) c_{\beta-1}(n) n^{\alpha+\beta-2} \right]
\nonumber\\
&& \quad \times 
\prod_{1 \leq i < j \leq N}
\Big\{(a_j-a_i)(b_j-b_i)\Big\}
\times \Big\{1+{\cal O}(|\a|, |\b|) \Big\}.
\nonumber
\end{eqnarray}
By (\ref{eqn:cpn}), we see
\begin{eqnarray}
&& \sum_{n \in \N} 
e^{-\pi^2 T n^2/\{2(\ell+r)^2\}}
c_{\alpha-1}(n) c_{\beta-1}(n) n^{\alpha+\beta-2} 
\nonumber\\
&=& \frac{(-1)^{(\alpha+\beta)/2}}
{4 (\alpha-1)! (\beta-1)!} 
\left( \frac{\pi}{\ell+r} \right)^{\alpha+\beta-2}
\sum_{n \in \N} 
e^{-\pi^2 T n^2/\{2(\ell+r)^2\}}
n^{\alpha+\beta-2}
\nonumber\\
&& \times
\left\{
\exp \left( \frac{2 \pi \sqrt{-1} \ell}{\ell+r} n \right)
+(-1)^{\alpha}+(-1)^{\beta}
+(-1)^{\alpha+\beta}
\exp \left( -\frac{2 \pi \sqrt{-1} \ell}{\ell+r} n \right) \right\}
\nonumber\\
&=&  \frac{(-1)^{(\alpha+\beta)/2}}
{4 (\alpha-1)! (\beta-1)!} 
\left( \frac{\pi}{\ell+r} \right)^{\alpha+\beta-2}
\nonumber\\
&& \times
\left\{ \sum_{n \in \Z}
\exp \left(-\frac{\pi^2 T}{2(\ell+r)^2} n^2
+\frac{2 \pi \sqrt{-1} \ell}{\ell+r} n \right)
n^{\alpha+\beta-2} \right.
\nonumber\\
&& \qquad \qquad \qquad
\left.
+(-1)^{\beta}\sum_{n \in \Z}
\exp \left(-\frac{\pi^2 T}{2(\ell+r)^2} n^2 \right)
n^{\alpha+\beta-2} \right\}.
\nonumber
\end{eqnarray}
By the identities (\ref{eqn:Iden1}),
the above is written as
\begin{eqnarray}
&& \frac{(-1)^{\alpha+\beta-1} (\ell+r)}
{2^{(\alpha+\beta+1)/2} \pi^{1/2}
(\alpha-1)! (\beta-1)!}
T^{-(\alpha+\beta-1)/2}
\nonumber\\
&& \times
\left\{ \Theta_{\alpha+\beta-2}
\left( \frac{2(\ell+r)}{\sqrt{2T}},
\frac{2\ell}{\sqrt{2T}} \right)
+(-1)^{\beta}
\Theta_{\alpha+\beta-2}
\left( \frac{2(\ell+r)}{\sqrt{2T}}, 0 \right) \right\}.
\nonumber
\end{eqnarray}
Therefore, we have the following asymptotics
in $|\a| \to 0$ and $|\b| \to 0$,
\begin{eqnarray}
\label{eqn:asym4}
&& f_N^{{\rm A}(-\ell, r)}(T, \b|\a)
\\
&=& \frac{(-1)^{N}}
{2^{N^2/2} \pi^{N/2} \prod_{k=1}^{N} \Gamma(k)^2}
T^{-N^2/2}
\nonumber\\
&& \times
\det_{1 \leq \alpha, \beta \leq N}
\left[ \Theta_{\alpha+\beta-2}
\left( \frac{2(\ell+r)}{\sqrt{2T}},
\frac{2\ell}{\sqrt{2T}} \right)
+(-1)^{\beta}
\Theta_{\alpha+\beta-2}
\left( \frac{2(\ell+r)}{\sqrt{2T}}, 0 \right)
\right]
\nonumber\\
&& \qquad \times
\prod_{1 \leq i < j \leq N}
\Big\{ (a_j-a_i)(b_j-b_i) \Big\} \times
\Big\{1+{\cal O}(|\a|, |\b|)\Big\}.
\nonumber
\end{eqnarray}

By the equalities (\ref{eqn:PEE}) and (\ref{eqn:equal1}), 
and using the asymptotics formulas (\ref{eqn:asym1})
and (\ref{eqn:asym4}),
we can take the limit
$$
\lim_{|\a|, |\b| \to 0}
\P \Big( -\ell < L^{\a, \b}(N,T),
R^{\a, \b}(N,T) < r \Big)
= \lim_{|\a|, |\b| \to 0}
\frac{f_N^{{\rm A}(-\ell, r)}(T, \b|\a)}
{f_N^{\rm A}(T, \b|\a)},
$$
and then the formula (\ref{eqn:thm2_1_1}) is obtained.

By the equalities (\ref{eqn:equal2}) and
using the asymptotics formulas (\ref{eqn:asym2})
and (\ref{eqn:asym3}),
we have
\begin{eqnarray}
\label{eqn:rep1}
&& \P \Big( -\ell < L^{\0, \R}(N,T),
R^{\0, \R}(N,T) < r \Big)
\\
&& = \lim_{|\a| \to 0}
\frac{\displaystyle{ \int_{\W_N^{\rm A}(-\ell, r)} d \b \,
f_N^{{\rm A}(-\ell, r)}(T, \b|\a)}}
{\cN_N^{\rm A}(T, \a)}
\nonumber\\
&& = \frac{(-1)^{N(N+1)/2}}
{2^{N(N+1)/4} \prod_{k=1}^{N} \Gamma(k/2)} T^{-N/2}
\nonumber\\
&& \quad \times \int_{\W_N^{\rm A}(-\ell,r)} d \b \,
\det_{1 \leq i, j \leq N} \left[
\Theta_{i-1} \left( \frac{2(\ell+r)}{\sqrt{2T}},
\frac{2\ell}{\sqrt{2T}}
+\frac{b_{j}}{\sqrt{2T}} \right)
\right.
\nonumber\\
&& \qquad \qquad \qquad \left.
+(-1)^{i}
\Theta_{i-1} \left( \frac{2(\ell+r)}{\sqrt{2T}},
\frac{b_{j}}{\sqrt{2T}} \right) \right]
\nonumber\\
&& = \frac{(-1)^{N(N+1)/2}}
{2^{N(N-1)/4} \prod_{k=1}^{N} \Gamma(k/2)} 
\nonumber\\
&& \quad \times 
\int_{\W_N^{\rm A}(-\ell/\sqrt{2T}, r/\sqrt{2T})} d \x \,
\det_{1 \leq i, j \leq N} \left[
\Theta_{i-1} \left( \frac{2(\ell+r)}{\sqrt{2T}},
\frac{2\ell}{\sqrt{2T}}
+x_j \right)
\right.
\nonumber\\
&& \qquad \qquad \qquad \left.
+(-1)^{i}
\Theta_{i-1} \left( \frac{2(\ell+r)}{\sqrt{2T}},
x_j \right) \right].
\nonumber
\end{eqnarray}

If we apply the following lemma,
which is known as the de Bruijn identity \cite{deBr55},
the formula (\ref{eqn:thm2_1_2}) is obtained.

\begin{lem}
\label{thm:lemma_deBruijn}
Let $z_i(x), 1 \leq i \leq \widehat{N}$ be an
integrable piecewise continuous function
on a region $\Lambda \subset \R$.
Let 
$
\W_N^{\rm A}(\Lambda)=
\Big\{ \x \in \W_N^{\rm A} :
x_i \in \Lambda, 1 \leq i \leq N \Big\}.
$
Then
\begin{equation}
\int_{\W_N^{\rm A}(\Lambda)} d \x \,
\det_{1 \leq i, j \leq N} 
\Big[ z_i(x_j) \Big]
= \Pf_{1 \leq i, j \leq \widehat{N}}
\Big[ Z_{ij} \Big],
\label{eqn:deBruijn}
\end{equation}
where
\begin{eqnarray}
&& I_i =
\int_{\Lambda} z_i(x) dx,
\nonumber\\
&& I_{ij}=
\int_{(x_1, x_2) \in \Lambda^2: x_1 < x_2} 
\det \left[ \begin{array}{ll}
z_i(x_1) & z_i(x_2) \cr
z_j(x_1) & z_j(x_2) 
\end{array} \right] dx_1 dx_2,
\nonumber
\end{eqnarray}
and
\begin{equation}
Z_{ij}
= \left\{ \begin{array}{rl}
I_{ij}, \quad & \mbox{if} \quad 1 \leq i, j \leq N,
\cr
I_{i}, \quad & \mbox{if} \quad 
1 \leq i \leq N, j=N+1, 
\cr
-I_{j}, \quad & \mbox{if} \quad
i=N+1, 1 \leq j \leq N,
\cr
0, \quad & \mbox{if} \quad
i=j= N+1.
\end{array} \right.
\nonumber
\end{equation}
\end{lem}

The proof of Theorem \ref{thm:Thm2_1} was then completed.

\SSC{Proof of Theorem \ref{thm:Thm2_2}}

For $h >0$, consider the restricted region in $\W_N^{\rm C}$,
$$
\W_N^{\rm C}(h)= \Big\{
\x=(x_1, \dots, x_N) \in \W_N^{\rm C} :
x_N < h \Big\}.
$$
The transition probability density of a one-dimensional
Brownian motion in an interval $(0, h)$,
where two absorbing walls are put at $x=0$ and $x=h >0$,
is given by
$$
p_{\rm abs}^{(0,h)}(t, y|x)
=\frac{2}{h} \sum_{n \in \N} e^{-\pi^2 t n^2/(2h^2)}
\sin \left( \frac{\pi x n}{h} \right)
\sin \left( \frac{\pi y n}{h} \right),
$$
$t >0, x,y \in \W_N^{\rm C}(h)$.
For $\a, \b \in \W_N^{\rm C}(h)$, let
$$
f_{N}^{{\rm C}(h)}(T, \b|\a)
=\det_{1 \leq i, j \leq N} \Big[
p_{\rm abs}^{(0,h)}(T, b_i|a_j) \Big].
$$
By the Karlin-McGregor formula, it gives the
probability density that the absorbing Brownian motion
in $\W_N^{\rm C}(h)$ starting from $\a \in \W_N^{\rm C}(h)$
at time 0 arrives at $\b \in \W_N^{\rm C}(h)$ at time $T$.

Then we have the equalities, for $h >0$
\begin{eqnarray}
\label{eqn:equal3}
&& \P \Big( H^{\a,\b}(N, T) < h \Big)
= \frac{f_{N}^{{\rm C}(h)}(T, \b|\a)}
{f_N^{\rm C}(T, \b|\a)},
\quad \a, \b \in \W_N^{\rm C}(h),
\\
\label{eqn:equal4}
&& \P \Big( H^{\a, \R_+}(N, T) < h \Big)
=\frac{\displaystyle{\int_{\W_N^{\rm C}(h)} d \b \,
f_N^{{\rm C}(h)}(T, \b|\a)}}
{\cN_N^{\rm C}(T, \a)},
\quad \a \in \W_N^{\rm C}(h).
\end{eqnarray}

The following asymptotics in $|\a| \to 0$ were reported in 
\cite{KT03a,KT04,KIK08b};
as $|\a| \to 0$,
\begin{eqnarray}
\label{eqn:asymC1}
&& f_N^{\rm C}(T, \b|\a)
= \frac{2^{N/2}}{\pi^{N/2} \prod_{m=1}^{N} \Gamma(2m)}
T^{-N(2N+1)/2} e^{-|\b|^2/(2T)}
\\
&& \quad \times
\prod_{1 \leq i<j \leq N}
\Big\{ (a_j^2-a_i^2)(b_j^2-b_i^2) \Big\}
\prod_{k=1}^{N} \Big\{ a_k b_k \Big\}
\times \Big\{1+{\cal O}(|\a|) \Big\}, 
\nonumber\\
\label{eqn:asymC2}
&& \cN_N^{\rm C}(T, \b|\a)
= \left( \frac{2}{\pi} \right)^{N/2}
\prod_{m=1}^{N} \frac{\Gamma(m)}{\Gamma(2m)} T^{-N^2/2}
\\
&& \quad \quad \times
\prod_{1 \leq i<j \leq N}
(a_j^2-a_i^2)
\prod_{k=1}^{N} a_k 
\times \Big\{1+{\cal O}(|\a|) \Big\}, 
\nonumber
\end{eqnarray}

By multilinearity of determinants
and by performing the series expansion
of a sine function in a determinant,
\begin{eqnarray}
\label{eqn:fNCh1}
&& f_N^{{\rm C}(h)}(T, \b|\a)
\\
&& \qquad 
= \frac{1}{N!} \left( \frac{2}{h} \right)^N
\sum_{\n \in \N^{N}} e^{-\pi^2 T |\n|^2/(2h^2)}
\nonumber\\
&& \qquad \qquad \times
\det_{1 \leq i, j \leq N}
\left[ \sin \left( \frac{\pi a_i}{h} n_j \right) \right]
\det_{1 \leq \alpha, \beta \leq N}
\left[ \sin \left( \frac{\pi b_{\alpha}}{h} n_{\beta} \right) \right]
\nonumber\\
&& \qquad = \frac{(-1)^{N(N-1)/2} 2^{N} \pi^{N^2}}
{h^{N(N+1)} \prod_{m=1}^{N} \Gamma(2m)}
\nonumber\\
&& \qquad \qquad \times
\det_{1 \leq \alpha, \beta \leq N}
\left[ \sum_{n \in \N}
e^{-\pi^2 T n^2/(2h^2)} n^{2\alpha-1}
\sin \left( \frac{\pi b_{\beta} n}{h} \right) \right]
\nonumber\\
&& \qquad \qquad \times
\prod_{1 \leq i<j \leq N}(a_j^2-a_i^2)
\prod_{k=1}^{N} a_k \times
\Big\{1+{\cal O}(|\a|) \Big\}.
\nonumber
\end{eqnarray}
Here we note that
\begin{eqnarray}
&& \sum_{n \in \N}
e^{-\pi^2 T n^2/(2h^2)} n^{2\alpha-1}
\sin \left( \frac{\pi b_{\beta} n}{h} \right)
\nonumber\\
&& \qquad 
= \frac{1}{2 \sqrt{-1}}
\sum_{n \in \Z} n^{2 \alpha-1}
\exp \left(
-\frac{\pi^2 T}{2 h^2} n^2
+\frac{\pi \sqrt{-1} b_{\beta}}{h} n \right)
\nonumber\\
&& \qquad 
= \frac{(-1)^{\alpha-1} h^{2 \alpha}}
{2^{\alpha} \pi^{(4 \alpha-1)/2}}
T^{-\alpha} 
\Theta_{2 \alpha-1} \left( \frac{2h}{\sqrt{2T}},
\frac{b_{\beta}}{\sqrt{2T}} \right),
\nonumber
\end{eqnarray}
where we have used the identity (\ref{eqn:Iden1}).
Then (\ref{eqn:fNCh1}) is written as
\begin{eqnarray}
\label{eqn:asymC3}
&& f_N^{{\rm C}(h)}(T, \b|\a)
\\
&& = \frac{T^{-N(N+1)/2}}
{2^{N(N-1)/2} \pi^{N/2} \prod_{m=1}^{N} \Gamma(2m)}
\det_{1 \leq \alpha, \beta \leq N}
\left[ \Theta_{2 \alpha-1} 
\left( \frac{2 h}{\sqrt{2T}},
\frac{b_{\beta}}{\sqrt{2T}} \right) \right]
\nonumber\\
&& \quad
\times \prod_{1 \leq i < j \leq N}
(a_j^2-a_i^2) \prod_{k=1}^{N} a_k
\times \Big\{1+{\cal O}(|\a|) \Big\}.
\nonumber
\end{eqnarray}

On the other hands, in $|\b| \to 0$,
(\ref{eqn:fNCh1}) has the asymptotics
\begin{eqnarray}
&& f_N^{{\rm C}(h)}(T, \b|\a)
\nonumber\\
&& = \frac{1}{N!} \left( \frac{2}{h}\right)^N
\frac{1}{\prod_{i=1}^{N} \Gamma(2m)^2}
\sum_{\n \in \N^N} e^{-\pi^2 T |\n|^2/(2h^2)}
\nonumber\\
&& \quad \times 
\det_{1 \leq \alpha, \beta \leq N}
\Big[ n_{\alpha}^{2\beta-1} \Big]
\det_{1 \leq \gamma, \delta \leq N}
\Big[ n_{\gamma}^{2\delta-1} \Big]
\nonumber\\
&& \quad \times
\prod_{1 \leq i<j \leq N} 
\Big\{(a_j^2-a_i^2)(b_j^2-b_i^2) \Big\}
\prod_{k=1}^N \Big\{ a_k b_k \Big\}
\times \Big\{ 1+{\cal O}(|\a|, |\b|) \Big\}
\nonumber\\
&& = \frac{2^N \pi^{2N^2}}
{h^{N(2N+1)} \prod_{m=1}^{N} \Gamma(2m)^2}
\det_{1 \leq \alpha, \beta \leq N}
\left[ \sum_{n \in \N} n^{2\alpha+2 \beta-2}
e^{-\pi^2 T n^2/(2h^2)} \right]
\nonumber\\
&& \quad \times
\prod_{1 \leq i<j \leq N} 
\Big\{(a_j^2-a_i^2)(b_j^2-b_i^2) \Big\}
\prod_{k=1}^N \Big\{ a_k b_k \Big\}
\times \Big\{ 1+{\cal O}(|\a|, |\b|) \Big\}.
\nonumber
\end{eqnarray}
By using the identities (\ref{eqn:Iden1}),
the above is written as follows,
\begin{eqnarray}
\label{eqn:asymC4}
&& f_N^{{\rm C}(h)}(T, \b|\a)
\\
&& = \frac{(-1)^{N}}
{2^{N(2N-1)/2} \pi^{N/2}
\prod_{m=1}^{N} \Gamma(2m)^2 }
T^{-N(2N+1)/2}
\nonumber\\
&& \quad \times
\det_{1 \leq \alpha, \beta \leq N}
\left[ \Theta_{2(\alpha+\beta-1) }
\left( \frac{2h}{\sqrt{2T}}, 0 \right) \right]
\nonumber\\
&& \quad \times
\prod_{1 \leq i<j \leq N} 
\Big\{(a_j^2-a_i^2)(b_j^2-b_i^2) \Big\}
\prod_{k=1}^N \Big\{ a_k b_k \Big\}
\times \Big\{ 1+{\cal O}(|\a|, |\b|) \Big\}.
\nonumber
\end{eqnarray}

Using the asymptotics formulas (\ref{eqn:asymC1}) and
(\ref{eqn:asymC4}), we can take the
limit $|\a|, |\b| \to 0$ in
(\ref{eqn:equal3}). Then
(\ref{eqn:thm2_2_1}) is obtained.

By the equality (\ref{eqn:equal4}), and using
the asymptotics (\ref{eqn:asymC2}) and (\ref{eqn:asymC3}),
we have
\begin{eqnarray}
&& \P \Big( H^{\0, \R_+}(N,T) < h \Big)
\nonumber\\
&& = \lim_{|\a| \to 0}
\frac{\displaystyle{ \int_{\W_N^{\rm C}(h)} d \b \,
f_N^{{\rm C}(h)}(T, \b|\a)}}
{\cN_N^{\rm C}(T, \a)}
\nonumber\\
&& = \frac{T^{-N/2}}{ 2^{N^2/2}
\prod_{k=1}^{N} \Gamma(k)}
\int_{\W_N^{\rm C}(h)} d \b \,
\det_{1 \leq i, j \leq N} \left[
\Theta_{2i-1} \left( \frac{2h}{\sqrt{2T}}, 
\frac{b_j}{\sqrt{2T}} \right) \right].
\nonumber\\
&& = \frac{1}{ 2^{N(N-1)/2} 
\prod_{k=1}^{N} \Gamma(k)}
\int_{\W_N^{\rm C}(h/\sqrt{2T})} d \x \,
\det_{1 \leq i, j \leq N} \left[
\Theta_{2i-1} \left( \frac{2h}{\sqrt{2T}}, x_j \right) \right].
\nonumber
\end{eqnarray}

Then, we use the de Bruijn identity
given in Lemma \ref{thm:lemma_deBruijn}
and the formula (\ref{eqn:thm2_2_2}) is obtained.
The proof of Theorem \ref{thm:Thm2_2} was then completed.

\vskip 1cm

\begin{small}
\noindent{\it Acknowledgements.} \quad
The present authors would like to thank Professor N. Minami
for giving them an opportunity for publishing the
present manuscript.
They also thank Dr.N. Kobayashi
for giving them the figures, which illustrate
samples of paths of noncolliding diffusion processes
discussed in the present paper.
M.K. is supported in part by
the Grant-in-Aid for Scientific Research (C)
(No.21540397) of Japan Society for
the Promotion of Science.


\end{small}

\begin{thebibliography}{99}
%
%
\bibitem{AZ97}
Altland, A. and Zirnbauer, M. R.,
Nonstandard symmetry classes 
in mesoscopic normal-superconducting hybrid structure,
\textit{Phys. Rev. B}, 
\textbf{55} (1997), 1142--1161.

\bibitem{AAR99}
Andrews, G. E., Askey, R. and Roy, R.,
\textit{Special Functions},
Cambridge University Press, Cambridge, 1999.

\bibitem{BPY01}
Biane, P., Pitman, J. and Yor, M.,
Probability laws related to the Jacobi theta and
Riemann zeta functions, and Brownian excursions,
\textit{Bull. Amer. Math. Soc.},
\textbf{38} (2001), 435--465.

\bibitem{BFPSW09}
Borodin, A., Ferrari, P.L., Pr\"ahofer, M.,
Sasamoto, T. and Warren, J.,
Maximum of Dyson Brownian motion and
non-colliding systems with a boundary,
\textit{Elect. Commun. in Probab.}
\textbf{14} (2009), 486-494.

\bibitem{deBr55}
de Bruijn, N. G.,
On some multiple integrals involving determinants,
\textit{J. Indian Math. Soc.} 
\textbf{19} (1955), 133--151.

\bibitem{Feierl07}
Feierl, T.,
The height of watermelons with wall (extended abstract),
\textit{Proceedings of the 
2007 Conference on Analysis of Algorithms},
Discrete Mathematics and Theoretical Computer Science
Proceedings, (2007).

\bibitem{Feierl09}
Feierl, T.,
The height and range of watermelons without wall
(extended abstract),
\textit{Proceedings of the International Workshop 
on Combinatorial Algorithms},
edited by J. Fiala, J. Kratochv\'il,
and M. Miller, 
\textit{Lecture Notes in Computer Science},
Springer-Verlag, Berlin, 
vol. 5874 (2009), 242--253.

\bibitem{Ful07}
Fulmek, M.,
Asymptotics of the average height of 2-watermelons
with a wall,
\textit{Electron. J. Comb.},
\textbf{14} (2007), {\#}R64/1--20.

\bibitem{FH91}
Fulton, W. and Harris, J.,
\textit{Representation Theory, A First Course},
Springer, New York, 1991.

\bibitem{KM59}
Karlin, S. and McGregor, J.,
Coincidence probabilities,
\textit{Pacific J. Math.}, 
\textbf{9} (1959), 1141--1164.

\bibitem{KIK08a}
Katori, M, Izumi, M. and Kobayashi, N.,
Two Bessel bridges conditioned never to collide,
double Dirichlet series, and Jacobi theta function,
\textit{J. Stat. Phys.},
\textbf{131} (2008), 1067--1083.

\bibitem{KT02}
Katori, M. and Tanemura, H.,
Scaling limit of vicious walks and two-matrix model,
\textit{Phys. Rev.} E 
\textbf{66} (2002), 011105/1--12.

\bibitem{KT03a}
Katori, M. and Tanemura, H.,
Functional central limit theorems for vicious walkers,
\textit{Stoch. Stoch. Rep.}, 
\textbf{75} (2003), 369--390.

\bibitem{KT03b}
Katori, M. and Tanemura, H.,
Noncolliding Brownian motions and Harish-Chandra formula,
\textit{Elect. Comm. in Probab.}, 
\textbf{8} (2003), 112--121.

\bibitem{KT04}
Katori, M. and Tanemura, H.,
Symmetry of matrix-valued stochastic processes and
noncolliding diffusion particle systems,
\textit{J. Math. Phys.}, 
\textbf{45} (2004), 3058--3085.

\bibitem{KT07a}
Katori, M. and Tanemura, H.,
Infinite systems of non-colliding generalized meanders
and Riemann--Liouville differintegrals,
\textit{Probab. Theory Relat. Fields}, 
\textbf{\bf 138} (2007), 113--156.

\bibitem{KT07b}
Katori, M. and Tanemura, H.,
Noncolliding Brownian motion and determinantal processes,
\textit{J. Stat. Phys.}, 
\textbf{129} (2007), 1233--1277.

\bibitem{KTNK03}
Katori, M.,Tanemura, H., Nagao, T. and Komatsuda, N.,
Vicious walk with a wall, noncolliding meanders,
chiral and Bogoliubov-de Gennes random matrices,
\textit{Phys. Rev.} E  
\textbf{\bf 68} (2003), 021112/1--16.

\bibitem{KIK08b}
Kobayashi, N., Izumi, M. and Katori, M.,
Maximum distributions of bridges of
noncolliding Brownian paths,
\textit{Phys. Rev.} E,
\textbf{78} (2008), 051102/1--15.

\bibitem{Mehta04}
Mehta, M. L.,
{\it Random Matrices}, 3rd ed.,
Elsevier Academic Press, London, 2004.

\bibitem{PM83}
Pandey, A. and Mehta, M. L.,
Gaussian ensembles of random Hermitian intermediate between
orthogonal and unitary ones,
\textit{Commun. Math. Phys.} 
\textbf{87}, (1983) 449--468.

\bibitem{RY98}
Revuz, D. and Yor, M.,
\textit{Continuous Martingales and Brownian Motion}, 3rd ed.,
Springer, New York, 1998.

\bibitem{SMCRF08}
Schehr, G., Majumdar, S. N., Comtet, A. and Randon-Furling, J.,
Exact distribution of the maximum height of
$p$ vicious walkers,
\textit{Phys. Rev. Lett.}
\textbf{\bf 101}, (2008) 150601/1--4.

\bibitem{Stem90}
Stembridge, J. R.,
Nonintersecting paths, pfaffians, and the plane partitions,
\textit{Adv. in Math.} 
\textbf{83} (1990), 96--131.

\bibitem{TW07}
Tracy, C. A. and Widom, H.,
Nonintersecting Brownian excursions,
\textit{Ann. Appl. Probab.}
\textbf{17} (2007), 953--979.
\end{thebibliography}
\end{document}